\documentclass [12 pt]{amsart}
\usepackage{amssymb,latexsym}
\usepackage{xypic}
\theoremstyle{plain}
\newtheorem{Thm}{Theorem}[section]
\newtheorem{Cor}{Corollary}[section]

\newtheorem{theorem}{Theorem}

\newtheorem{lemma}[theorem]{Lemma}

\numberwithin{equation}{section}
\theoremstyle{definition}
\newtheorem{Def}{Definition}[section]
\theoremstyle{remark}

\begin{document}

\title[On Fuzzification of $n$-Lie Algebras]{On Fuzzification of $n$-Lie Algebras
}

\author{Shadi Shaqaqha}

\address[Shadi Shaqaqha]{Department of Mathematics, Yarmouk University, Irbid, Jordan
}
\email[Shadi Shaqaqha]{shadi.s@yu.edu.jo}



\begin{abstract}
The aim of this paper is to introduce the notion of intuitionistic fuzzy Lie subalgebras and intutionistic fuzzy Lie ideals of
$n$-Lie algebras. It is a generalization of intuitionistic fuzzy Lie algebras. Then, we investigate some of characteristics of intutionistic fuzzy Lie ideals (resp. subalgebras) of $n$-Lie algeras. Finally, we define the image and the
inverse image of intuitionistic fuzzy Lie subalgebra under $n$-Lie algebra homomorphism. The properties of intuitionistic fuzzy $n$-Lie
subalgebras and intuitionistic fuzzy Lie ideals under homomorphisms of $n$-Lie algebras are studied. Finally, we define the intuitionistic fuzzy quotient $n$-Lie algebra by an
intuitionistic fuzzy ideal of n-Lie algebra and prove that it is a $n$-Lie algebra.
\end{abstract}

\maketitle

\section{Introduction}
The concept of $n$-Lie algebras was introduced by Filippov \cite{Filippov} in 1985. If $n=2$, then we get Lie algebra structure which was introduced by Sophus Lie (1842-1899) while he was attempting to classify certain smooth subgroups of general linear groups that are now called Lie groups. The case where $n=3$ was initially appeared in Nanmbu's work \cite{Nambu} while he was attempting to describe simultaneous classical dynamics of three particles. Takhtajan \cite{Takhtajan} investigated the geometrical and algebraic aspects of the generalized Nambu mechanics, and established the connection between the Nambu mechanics and Filippov’s theory
of $n$-Lie algebras. More applictions to the $n$-Lie algebras can be found in \cite{Alekseevsky, Gautheron, Ho, Marmo, Michor,  Nakanishi, Papadopoulos, Vinogradov}. So, our results are expected to be useful in various applications.\\
The notion of fuzzy sets was firstly introduced by Zadeh \cite{Zadeh}. The fuzzy set theory states that there are propositions with an infinite number of truth values, assuming two extreme values, $1$ (totally true), $0$ (totally false) and a continuum in between, that justify
the term fuzzy. Applications of this
theory can be found, for example, in artificial intelligence, computer science, control engineering, decision
theory, logic and management science.\\
 After introducing of fuzzy sets by Zadeh, many researches were conducted on the generalizations of this fundamental concept. Among these generalizations was the concept of intuitionistic fuzzy sets, which was introduced by Attanasov \cite{ATANASSOV} in 1986, is the most important and interesting one. The elements of the intuitionistic fuzzy sets are distinguished by an additional degree called the degree of uncertainty. There are numerous applications to this new concept include computer science, mathematics, medicine, chemistry, economics, astronomy etc.. Many mathematicians have involved in extending the concept of intuitionistic fuzzy sets to border of abstract algebra.  Biswas applied the concepts of intuitionistic fuzzy sets to the theory of groups and studied intuitionistic fuzzy subgroups in \cite{Biswas}. Also, Akram studied Lie algebra in intuitionistic fuzzy sets and obtained some results in \cite{Akram}.\\
The fuzzy Lie subalgebras and fuzzy Lie ideals are considered in \cite{Kim} by Kim and
Lee, and in \cite{Yehia1,Yehia2} by Yehia. They established the analogues of most of the fundamental ground results
involving Lie algebras in the fuzzy setting. The study of fuzzy subalgebras (resp. ideals) of $n$-Lie algebras was initiated by B. Davvaz and WA. Dudek \cite{Davvaz}. Recently, the complex (intuitionistic) fuzzy Lie algebras is studied in \cite{shadi, shadi2} as
a generalization of (intutionistic) fuzzy Lie algebras.\\
In this paper we describe intuitionistic fuzzy n-Lie algebras. We will introduce n-Lie
algebras into intuitionistic fuzzy set. Our work will generalize the theory of (intu-
itionstic) fuzzy Lie algebras (\cite{Akram, Akram1, Akram2, Akram3, shadi2, Yehia1, Yehia2}).
\section{Preliminaries}
$n$-Lie algebras were originally introduced by Filippov \cite{Filippov} in 1985. They generalize Lie algebras. In this article the ground field $F$ is arbitrary.
\begin{Def}\label{L1}{
Let $n\in \mathbb{N}$, $n\geq 2$. An $n$-Lie algebra is a pair $(L, [~])$ where $L$ is a vector space and 
$$[~]: L^n\rightarrow L;~(x_1, \ldots, x_n)\mapsto [x_1, \ldots, x_n]$$
is an $n$-linear map, called $n$-Lie bracket, that satisfies the following identities for all $\sigma$ in the symmetric group $S_n$ and $x_1, \ldots, x_n, y_2, \ldots, y_n\in L$:
\begin{itemize}
\item[(i)] Skew symmetry:
$$[x_{\sigma(1)}, \ldots, x_{\sigma(n)}]=\mathrm{sign}(\sigma)[x_1, \ldots, x_n].$$
\item[(ii)] The generalized Jacobi identity (called also the Filippov
identity):
$$[[x_1, \ldots, x_n],y_2, \ldots, y_n]=\sum_{i=1}^n[x_1, \ldots, x_{i-1}, [x_i, y_2, \ldots, y_n], x_i, \ldots, x_n].$$
\end{itemize}
}
\end{Def}
Subalgebras of $n$-Lie algebras and homomorphisms (or isomorphisms) between $n$-Lie algebras are defined as usual (\cite{Filippov}).\\
Throughout this paper, $L$ is a $n$-Lie algebra over field $F$.\\
The concept of intuitionistic fuzzy set  was introduced by Atanassov \cite{ATANASSOV}, where he added a new component (which determines the degree of non-membership) in the definition of fuzzy set (FS) that was given by Zadeh. Let $X$ be a non-empty set, and let $A=(\mu_A, \lambda_A)=\{(x, \mu_A(x), \lambda_A(x))~:~x\in X\}$  where $\mu_A: X\rightarrow [0, 1]$ and $\lambda_A: X\rightarrow [0, 1]$ be mappings such that $\mu_A(x)+\lambda_A(x)\leq 1$. Then $A$ is called an intuitionistic fuzzy set. In this case the mappings $\mu_A$ and $\lambda_A$ denote the degree of membership and the degree of non-membership to $A$ respectively, for each element $x\in X$. The value $\pi_A(x)=1-\mu_A(x)-\lambda_A(x)$ is called uncertainty or intuitionistic index of the element $x\in X$ to the intuitionistic fuzzy set $A$. It is obvious that each fuzzy set $A=\mu_A=\{(x, \mu_A(x))~:~x\in X\}$ can be represented as an intuitionistic fuzzy set where $A=\{(x, \mu_A(x), 1-\mu_A(x))~:~x\in X\}$.
\begin{Def}\label{IFSR}
Let $A=(\mu_A, \lambda_A)$ and $B=(\mu_B, \lambda_B)$ be two intuitionistic fuzzy sets of a set $X$. Then
\begin{itemize}
    \item[(i)] The complement of $A$ is $\Bar{A}= (\lambda_A, \mu_A)$,
    \item[(ii)] $A\subseteq B$ if and only if $\mu_A(x)\leq \mu_B(x)$ and $\lambda_A(x)\geq \lambda_B(x)$ for all $x\in X$,
    \item[(iii)] the intersection of $A$ and $B$ is $A\cap B = \{(x, \mathrm{min}\{\mu_A(x), \mu_B(x)\}, \mathrm{max}\{\lambda_A(x), \lambda_B(x)\}~|~x\in X$,
    \item[(iv)] the union of $A$ and $B$ is $A\cup B= \{(x, \mathrm{max}\{\mu_A(x), \mu_B(x)\}, \mathrm{min}\{\lambda_A(x), \lambda_B(x)\}~|~x\in X$,
    \item[(v)] $\Box A=\{(x,~\mu_A(x),~\mu_A^\complement(x)):~x\in X\}$ where $\mu_A^\complement(x)=(1-\mu_A(x))$ for all $x\in X$,
    \item [(vi)] $\diamondsuit A=\{(x,~\lambda_A^\complement (x),~\lambda_A(x)):~x\in X\}$ where $\lambda_A^\complement(x)=(1-\lambda_A(x))$ for all $x\in X$.
\end{itemize}
\end{Def}
Let $L$ be a $n$-Lie algebra over $F$. A fuzzy set $A=\mu_A=\{(x,~\mu_A(x)):x\in L\}$  on $L$ is a fuzzy Lie subalgebra if the following conditions are satisfied:
\begin{itemize} 
\item[(i)] $\mu_{A}(x+y) \geq \mathrm{min}\{\mu_{A}(x), \mu_{A}(y)\}$ for all $x, y\in L$,
\item[(ii)] $\mu_{A}(\alpha x) \geq \mu_{A}(x)$ for all $x\in L$ and $\alpha \in F$,
\item[(iii)] $\mu_{A}([x_1,~x_2, \ldots, x_n]) \geq \mathrm{min}\{\mu_{A}(x_1), \mu_{A}(x_2),  \cdots, \mu_A(x_n)\}$ for all $ x_1,~x_2, \ldots, ~x_n \in L$.
\end{itemize}
It is called a fuzzy Lie ideal if the condition
$$\mu_{A}([x_1,~x_2, \ldots, x_n]) \geq \mathrm{min}\{\mu_{A}(x_1), \mu_{A}(x_2),  \cdots, \mu_A(x_n)\}$$
is replaced by
$$\mu_{A}([x_1,~x_2, \ldots, x_n]) \geq \mathrm{max}\{\mu_{A}(x_1), \mu_{A}(x_2),  \cdots, \mu_A(x_n)\}$$
(\cite{Davvaz}).
\section{Intuitionistic Fuzzy $n$-Lie Algebras}
For the sake of simplicity, we shall use the symbols $a\wedge b =\mathrm{min}\left\{a, b\right\}$ and $a\vee b =\mathrm{max}\left\{a, b\right\}$.
\begin{Def}\label{ifls}
	Let $A=(\mu_A,\lambda_A)=\{(x,~\mu_A(x),~\lambda_A(x)):x\in L\}$ be an intuitionistic fuzzy set of $L$. Then $A$ is called an intuitionistic fuzzy Lie subalgebra of $L$ if the following conditions are satisfied:
	\begin{itemize}
		\item [(i)] $\mu_{A}(x+y) \geq \mu_{A}(x)\wedge \mu_{A}(y)$ and $\lambda_{A}(x+y) \leq \lambda_{A}(x)\vee\lambda_{A}(y)$ for all $x, y\in L$,
		\item[(ii)] $\mu_{A}(\alpha x) \geq \mu_{A}(x)$ and $\lambda_{A}(\alpha x) \leq \lambda_{A}(x)$ for all $x\in L$ and $\alpha \in F$,
		\item[(iii)] $\mu_{A}([x_1,~x_2, \ldots,~x_n]) \geq \mu_{A}(x_1)\wedge\mu_{A}(x_2)\wedge \cdots \wedge \mu_A(x_n)$ and  $\lambda_{A}([x_1,~x_2, \ldots,~x_n]) \leq \lambda_{A}(x_1)\vee\lambda_{A}(x_2)\vee \cdots \vee \lambda_A(x_n)$ for all $x_1, x_2, \ldots, x_n\in L$.
	\end{itemize}
	An intuitionistic fuzzy set $A$ on $L$ is called an intuitionistic fuzzy Lie ideal if the conditions $(i)$ and $(ii)$ are satisfied together with the following addition condition:\\
	$(iii)^\prime \mu_{A}([x_1,~x_2, \ldots,~x_n]) \geq \mu_{A}(x_1)\vee\mu_{A}(x_2)\vee \cdots \vee \mu_A(x_n)$ and  $\lambda_{A}([x_1,~x_2, \ldots,~x_n]) \leq \lambda_{A}(x_1)\wedge\lambda_{A}(x_2)\wedge \cdots \wedge \lambda_A(x_n)$ for all $x_1, x_2, \ldots, x_n\in L$.
\end{Def}
In the special case that $n=2$, we obtain intuitionistic fuzzy Lie algebras (\cite{Akram}). When first two conditions hold, we say that $A$ is an intuitionistic fuzzy vector subspace of $L$. The second condition implies $\mu_A(0)\geq \mu_A(x)$, $\mu_A(-x)\geq \mu_A(x)$, $\lambda_A(0)\leq \lambda_A(x)$, and $\lambda_A(-x)\leq \lambda_A(x)$ for all $x\in L$. Also, every intuitionistic fuzzy Lie ideal of an $n$-Lie algebra is an intuitionistic fuzzy Lie subalgebra, but the converse is not necessary true ( see \cite[Example 3.1]{shadi}).\\

Let $\{A_i=(\mu_{A_i}, \lambda_{A_i})~|~i\in I\}$ be a collection of intuitionistic fuzzy sets on a nonempty set $X$. Then
$$\bigcap_{i\in I}A_i=\{(x,~\mu_{\bigcap_{i\in I}A_i}(x),~\lambda_{\bigcap_{i\in I}A_i} (x)):~x\in X \},$$
where $$\mu_{\bigcap_{i\in I}A_i}(x)=\inf_{i\in I}\{\mu_{A_i}(x)\}$$ and $$\lambda_{\bigcap_{i\in I}A_i} (x)=\sup_{i\in I}\{\lambda_{A_i}(x)\},$$
is an intuitionistic fuzzy set too (see \cite{shadi2}). We shall give the proof of the following theorem, established in \cite{shadi} to the case of intuitionistic fuzzy Lie algebras, which proves that the arbitrary intersection of intuitionistic fuzzy Lie subalgebras (resp. ideals) of an $n$-Lie algebra $L$ is an intuitionistic fuzzy Lie subalgebras (resp. ideal) of $L$ too. However it was proved in the case that $L$ is a Lie (super)algebras and where the family is finite (see \cite{chen2}).
\begin{Thm}\label{iflt}
	Let $\{A_i\}_{i\in I}$ $(A_i=(\mu_{A_i}, \lambda_{A_i}), i\in I)$ be a collection of intuitionistic fuzzy Lie subalgebras (resp. ideals) on $L$. Then $\bigcap_{i\in I}A_i$ is an intuitionistic fuzzy Lie subalgebra (resp. ideal) of $L$.
\end{Thm}
{\it Proof.~}Here we will prove the case of intuitionistic fuzzy Lie subalgera. For $x,~y\in L$ and $\alpha \in F$, we have
\begin{align*}
\mu_{\bigcap_{i\in I}A_i}(x+y)=&\inf_{i\in I}\{\mu_{A_i}(x+y)\}\\
\geq&\inf_{i\in I}\{\mu_{A_i}(x)\wedge\mu_{A_i}(y)\}\\
=&\inf_{i\in I}\{\mu_{A_i}(x)\}\wedge\inf_{i\in I}\{\mu_{A_i}(y)\}\\
=&\mu_{\bigcap_{i\in I}A_i}(x)\wedge\mu_{\bigcap_{i\in I}A_i}(y).
\end{align*}
Also,
\begin{align*}
\mu_{\bigcap_{i\in I}A_i}(\alpha x)=&\inf_{i\in I}\{\mu_{A_i}(\alpha x)\}\\
\geq&\inf_{i\in I}\{\mu_{A_i}(x)\}\\
=&\mu_{\bigcap_{i\in I}A_i}(x).
\end{align*}
Similarly, we can prove that $\lambda_{\bigcap_{i\in I}A_i}(x+y)\leq \lambda_{\bigcap_{i\in I}A_i}(x)\vee\lambda_{\bigcap_{i\in I}A_i}(y)$ and $\lambda_{\bigcap_{i\in I}A_i}(\alpha x)\leq \lambda_{\bigcap_{i\in I}A_i}(x)$. Next, if $x_1, \ldots, x_n\in L$, then
\begin{align*}
\mu_{\bigcap_{i\in I}A_i}([x_1,\ldots, x_n])=&\inf_{i\in I}\{\mu_{A_i}([x_1, \ldots, x_n])\}\\
\geq&\inf_{i\in I}\{\mu_{A_i}(x_1)\wedge \cdots \wedge \mu_{A_i}(x_n)\}\\
=&\inf_{i\in I}\{\mu_{A_i}(x)\}\wedge \cdots \wedge \inf_{i\in I}\{\mu_{A_i}(x_n)\}\\
=&\mu_{\bigcap_{i\in I}A_i}(x_1)\wedge \cdots \wedge \mu_{\bigcap_{i\in I}A_i}(x_n).
\end{align*}
In a similar way, one can show $\lambda_{\bigcap_{i\in I}A_i}([x_1, \ldots, x_n]) \leq \lambda_{\bigcap_{i\in I}A_i}(x)\vee \cdots \vee \lambda_{\bigcap_{i\in I}A_i}(x_n)$. Therefore, $\bigcap_{i\in I}A_i$ is an intuitionistic fuzzy Lie subalgebra. The proof of the case of intuitionistic fuzzy Lie ideal is same, so we omit it. \hfill$\Box$

Let $A=\{(x,~\mu_A(x),~\lambda_A(x)):~x\in X\}$ be an intuitionistic fuzzy set. For $s, t\in [0,~1] $ the set $A^{t}_{s}=\{x\in X:\mu_A(x)\geq s,~\lambda_A(x)\leq t\}$ is called the upper level subset of the intuitionistic fuzzy subset $A$. In particular if $t=1$, then we get the upper $s$-level cut $A_{s}^1=U(\mu_A;s)=\{x\in X~:~\mu_A(x)\geq s\}$. Also, if $s=0$, then we get the lower $t$-level cut $A_{0}^{t}= L(\lambda_A;t)=\{x\in X~:~\lambda_A(x)\leq t\}$.\\
The following two theorems will show relations between intuitionistic fuzzy Lie subalgebras of $L$ and Lie subalgebras of $L$. They are very similar to the case that
 suggested by Kondo and Dudek in \cite{Kondo}.
\begin{Thm}\label{IFL3}
	Let $A=(\mu_A,~\lambda_A)$ be an intuitionistic fuzzy set of an $n$-Lie algebra $L$.  Then A is an intuitionistic fuzzy Lie subalgebra of $L$ if and only if the non-empty set $A^{t}_{s}$ is Lie subalgebra for all $s, t\in [0,~1]$.
\end{Thm}
{\it Proof.~}Let  $A=\{(x, ~\mu_A(x),~ \lambda_A(x)):~x\in L\}$ be an intuitionistic fuzzy Lie subalgebra. For $x,~y\in A^{s}_{t}$ and $\gamma \in F$
\begin{itemize}
\item[(i)] $\mu_A(x+y)\geq \mu_A(x)\wedge \mu_A(y)\geq s$ and $\lambda_A(x+y)\leq \lambda_A(x)\vee \lambda_A(y)\leq t$,
\item[(ii)] $\mu_A(\gamma x)\geq \mu_A(x)\geq s$ and $\lambda_A(\gamma x)\leq \lambda_A(x)\leq t $.
\end{itemize}
Thus $x+y,~\gamma x\in A_s^t$. Also for $x_1, \ldots, x_n\in L$, we have $\mu_A([x_1,, \ldots, x_n])\geq\mu_A(x_1)\wedge \cdots \wedge \mu_A(x_n\geq s$ and $\lambda_A([x_1,\ldots, x_n])\leq \lambda_A(x_1)\vee \cdots \vee \lambda_A(x_n)\leq t$. That is $[x_1, \ldots, x_n]\in A^{t}_{s}$. Therefore $A^t_s$ is Lie subalgebra of $L$. Conversely, suppose that $A^s_t\neq \emptyset$ is a Lie subalgebra of $L$ for every $s, t\in [0,~1]$.  Let $x, y\in L$ and $\alpha\in F$. Fix $s_1=\mu_A(x)\wedge \mu_A(x)$ and $t_1=\lambda_A(x)\vee \lambda_A(y)$, so that $x, y\in A_{s_1}^{t_1}$. Since $A_{s_1}^{t_1}$ is a subspace of $L$, we have $x+y$ and $\alpha x$ are in $A_{s_1}^{t_1}$, and so $\mu_A(x+y)\geq s_1= \mu_A(x)\wedge \mu_A(y)$, $\lambda_A(x+y)\leq t_1= \lambda_A(x)\vee \lambda_A(y)$. Also for $x\in L$, set $s_1=\mu_A(x)$ and $t_1=\lambda_A(x)$. Then $x\in A_{s_1}^{t_1}$, and so $\gamma x\in A^{t_1}_{s_1}$. Hence $\mu_A(\gamma x)\geq \mu_A(x)$ and $\lambda_A(\gamma x)\leq \lambda_A(x)$. Finally let $x_1, x_2, \ldots , x_n\in L$. Fix $t= \mu_A(x_1)\wedge \cdots \wedge \mu_A(x_n)$ and $s= \lambda_A(x_1)\vee \cdots \vee \lambda_A(x_n)$. Thus $x_i\in A^t_s$ for all $i=1, \ldots, n$. Since $A^t_s$ is a subalgebra of $L$, we have $[x_1, \ldots, x_n]\in A^t_s$, so that $\mu_A([x_1, \ldots, x_n])\geq s = \mu_{A}(x_1)\wedge \cdots \wedge \mu_A(x_2)$ and $\lambda_A([x_1, \ldots, x_n])\leq t= \lambda_A(x_1)\vee \cdots \vee \lambda_A(x_n)$. The case of intuitionistic fuzzy ideals is almost same.
\hfill $\Box$
%Similarly, one can prove the following theorem.
%\begin{Thm}\label{CFL4}
%Let $\mu_A$ be a complex fuzzy subset of $L$. Then the following statements are equivalent:
%\begin{itemize}
%\item[(i)] $\mu_A$ is a complex fuzzy ideal of $L$,
%\item[(ii)] the upper level $U(\mu_A, t)$ is an ideal of $L$ for every $t\in \mathrm{Im}(\mu_A)$.
%\end{itemize}
%\end{Thm}

Let $A=\{(x,~\mu_A(x),~\lambda_A(x)):~x\in X\}$ be an intuitionistic fuzzy set. For $s, t\in [0,~1] $ the set $A^{t^<}_{s^>}=\{x\in X:\mu_A(x)> s,~\lambda_A(x)< t\}$ is called the strong upper level subset of the intuitionistic fuzzy subset $A$. 
\begin{Thm}\label{IFL4}
	Let $A=( \mu_A,~ \lambda_A)$ be an intuitionistic fuzzy subset of $L$. Then $A$ is an intuitionistic fuzzy Lie subalgebra of $L$ if and only if the non empty set $A^{t^<}_{s^>}$ is Lie subalgebra, for all $s, t\in [0,~1]$.
\end{Thm}
{\it Proof.~}The proof of the forward direction is almost identical to the proof in Theorem \ref{IFL3}. Conversely, suppose that the non-empty set $A^{t^<}_{s^>}$ is a Lie subalgebra for all $s, t\in [0,~1] $.  We need to show that the conditions of Definition \ref{ifls} are satisfied. Let $x, y\in L$. If $\mu_A(x)=0$ or $\mu_A(y)=0$, then it is clear that $\mu_A(x+y)\geq \mu_A(x)\wedge \mu_A(y)$, so we may assume that $\mu_A(x)\neq 0$ and $\mu_A(y)\neq 0$. Let $s_0$ be the largest number on the interval $[0, 1]$ such that $s_0<\mu_A(x)\wedge\mu_A(y)$ and there is no $a\in L$ satisfying $s_0< \mu_A(a)< \mu_A(x)\wedge\mu_A(y)$. Having $x, y\in A_{s_0^>}^{t^<}$ where $1>t>\lambda_A(x)\vee \lambda_A(y)$ (such $t$ exists because $\mu_A(x)$ and $\mu_A(y)$ are greater than $0$ in addition to $\mu_A(x)+ \lambda_A(x)$ and $\mu_A(y)+ \lambda_A(y)$ are less than or equals to $1$) implies that $x+y\in A^{t}_{s_0}$, and hence $\mu_A(x+y)> s_0$. Since there exist no $a\in L$ with $s_0< \mu_A(a)< \mu_A(x)\wedge \mu_A(y)$, it follows that $\mu_A(x+y)\geq\mu_A(x)\wedge \mu_A(y)$. If $\lambda_A(x)=1$ or $\lambda_A(y)= 1$, then it is obvious that $\lambda_A(x+y)\leq \lambda_A(x)\vee \lambda_A(y)$. Assume now that $\lambda_A(x)\neq 1\neq \lambda_A(y)$. Let $t_0$ be the smallest number on the interval $[0, 1]$ such that $t_0> \lambda_A(x)\vee \lambda_A(y)$ and there is no $a\in L$ with $t_0> \lambda_A(a)> \lambda_A(x)\vee \lambda_A(y)$. Hence $x, y\in A^{t_0^<}_{s^>}$ where $s< \mu_A(x)\wedge \mu_A(y)$, and so $x+y\in A^{t_0^<}_{s^>}$. Therefore $\lambda_A(x+y)< t_0$. Since there is no $a\in L$ such that $t_0> \lambda_A(a)> \lambda_A(x)\vee \lambda_A(y)$, it follows that $\lambda_A(x+y)\leq \lambda_A(x)\vee \lambda_A(y)$ as desired. In a similar way we can show that $\mu_A(\gamma x)\geq \mu_A(x)$ and $\lambda_A(\gamma x)\leq \lambda_A(x)$ for all $x\in L$ and $\gamma\in F$. Next let $x_1, \ldots, x_n\in L$.  Again we may assume that neither of $\mu_A(x_1), \ldots, \mu_A(x_n)$ is $0$. Let $s_1$ be the greatest number on the interval $[0, 1]$ such that $s_1<\mu_A(x)\wedge \cdots \wedge \mu_A(x_n)$ and there is no $a\in L$ such that $s_1< \mu_A(a)< \mu_A(x_1) \wedge \cdots \wedge \mu_A(x_n)$. Then $x_i\in A^{t^<}_{s_1^>}$, where $1>t> \lambda_A(x_1)\vee \cdots \vee \mu_A(x_n)$, for each $i=1, \ldots, n$. Since $A^{t^<}_{s_1^>}$ is subalgebra of $L$, we have $[x_1, \ldots , x_n]\in A^{t^<}_{s_1^>}$. Hence $\mu_A(x_1, \ldots, x_n)\geq \mu_A(x_1)\wedge \cdots \wedge \mu_A(x_n)$. In a similar fashion, we can prove that $\lambda_A([x_1, \ldots, x_n])\leq \lambda_A(x_1)\vee \cdots \vee \lambda_A(x_n)$ for all $x_1, \ldots , x_n\in L$. \hfill $\Box$

From the proofs of Theorem \ref{IFL3} and Theorem \ref{IFL4}, we can immediately obtain the following result. 
\begin{Cor}\label{IFL5}
Let $A=(\mu_A, \lambda_A)$ be an intuitionistic fuzzy subset of an $n$-Lie algebra $L$. The following statements are equivalent for every $s, t\in [0, 1]$:
\begin{itemize}
    \item[(i)] $A$ is an intuitionistic fuzzy Lie subalgebra (resp. ideal) of $L$,
    \item[(ii)] The nonempty subsets $A^t_s$ are $n$-Lie subalgebras (resp. ideals) of $L$,
    \item[(iii)] The nonempty subsets $A^{t^<}_{s}=\{x\in L~|~ \mu_A(x)\geq s ~\mathrm{and}~ \lambda_A(x)< t\}$ are $n$-Lie subalgebras (resp. ideals) of $L$,
    \item[(iv)] The nonempty subsets $A^t_{s^>}= \{x\in L~|~\mu_A(x)> s~\mathrm{and}~\lambda_A(x)\leq t\}$ are $n$-Lie subalgebras (resp. ideals) of $L$,
    \item[(v)] The The nonempty subsets $A^{t^<}_{s^>}$ are $n$-Lie subalgebras (resp. ideals) of $L$.
\end{itemize}
\end{Cor}
\section{Operations On Intuitionistic Fuzzy $n$-Lie ideals}
We omit the proofs for the following two results because they are straightforward.
\begin{Thm}\label{boxifnl1}
	Let $A=(\mu_A,~\lambda_A)$ be an intuitionistic fuzzy set in an $n$-Lie algebra $L$. Then $A$ is an intuitionistic fuzzy Lie subalgebra (resp. ideal) if and only if $\Box A$ and $\diamondsuit A$ are intuitionistic fuzzy Lie subalgebras (resp. ideals).
\end{Thm}
Let $A=( \mu_A,~ \lambda_A)$ and $B=( \mu_B,~ \lambda_B)$ be two  intuitionistic fuzzy subsets on a $n$-Lie algebra $L$. Let us introduce the following sum of $A$ and $B$, which was defined by Chen and Zhang \cite{chen2} in the case of intuitionistic fuzzy Lie superalgebras:
$$A+B=(\mu_{A+B},~\lambda_{A+B}),$$
where
$$\mu_{A+B}(x)=\sup_{x=a+b}\{\mu_A(a)\wedge\mu_B(b)\},$$ and $$\lambda_{A+B} (x)=\inf_{x=a+b}\{\lambda_A(a)\vee\lambda_B(b)\}.$$ 
If $A=( \mu_A,~ \lambda_A)$ and $B=( \mu_B,~ \lambda_B)$ are intuitionistic fuzzy sets on $L$, then $A+B$ is an intuitionistic fuzzy set of $L$. Indeed if $x\in L$ with $\mu_{A+B}(x)+\lambda_{A+B}(x)>1$, then $\sup_{x=a+b}\{\mu_A(a)\wedge \mu_B(b)\}+\inf_{x=a+b}\{\lambda_A(a)\vee\lambda_B(b)\}>1$, and so there exist $a_1,~b_1\in L$ such that $x=a_1+b_1$ with $\mu_A(a_1)\wedge \mu_B(b_1)+\inf_{x=a+b}\{\lambda_A(a)\vee\lambda_B(b)\}>1$. It follows that $\mu_A(a_1)\wedge \mu_B(b_1)+\lambda_A(a_1)\vee\lambda_B(b_1)>1$. Without loss of generality we may assume $\mu_A(a_1)\leq \mu_B(b_1)$. If $\lambda_A(a_1)\geq\lambda_B(b_1)$, then it is a clear contradiction. If $\lambda_A(a_1)\leq\lambda_B(b_1)$, then $\mu_A(a_1)+\lambda_B(b_1)>1$. Hence $\mu_A(a_1)+1-\mu_B(b_1)>1$, and so $\mu_A(a_1)>\mu_B(b_1)$. Contradiction. Therefore, $A+B$ is an intuitionistic fuzzy set of $L$.

It is well known that if $A$ and $B$ are ideals of $L$, then $A+B$ is an ideal of $L$ too. We are going to obtain an analogous result in the case of intuitionistic fuzzy Lie ideals. Similar result was obtained for complex fuzzy Lie subalgebra in \cite{shadi}, and for intuitionistic fuzzy Lie sub-superalgebras in \cite{chen2}.
\begin{Thm}\label{sumt}
	Let $A=( \mu_A,~ \lambda_A)$ and $B=( \mu_B,~ \lambda_B)$ be two intuitionistic fuzzy $n$-Lie ideals on $L$. Then $A+B$ is an intuitionistic fuzzy $n$-Lie ideal of $L$ too.
\end{Thm}
{\it Proof.~}We proceed as in the proof of corresponding result on complex fuzzy Lie algebras \cite{shadi2}. The only difference appears in the proof is to show that $\mu_A([x_1, \ldots, x_n])\geq \mu_A(x_1)\vee \cdots \vee \mu_A(x_n)$ and $\lambda_A([x_1, \ldots, x_n])\geq \lambda_A(x_1)\wedge \mu_A(x_n)$ for each $x_1, \ldots, x_n\in L$. Let $x_1, \ldots, x_n\in L$. Suppose that $\mu_{A+B}([x_1, \ldots, x_n])<\mu_{A+B}(x_1)\vee\cdots \vee \mu_{A+B}(x_n)$. Then there exists $t$ such that $\mu_{A+B}([x_1, \ldots, x_n])<t<\mu_{A+B}(x_1)\vee \ldots \vee \mu_{A+B}(x_n)$. Without loss of generality we may assume $\mu_{A+B}(x_1)=\mu_{A+B}(x_1)\vee \cdots \vee \mu_{A+B}(x_n)$. Thus $\mu_{A+B}([x_1, \ldots, x_n])<t<\sup_{x_1=a+b}\{\mu_A(a)\wedge\mu_B(b)\}$, and so there exist $a_1,~b_1\in L$ such that $\mu_{A+B}([x,~y])<t<\mu_A(a_1)\wedge\mu_B(b_1)$. Therefore $\mu_A(a_1)>t$ and $\mu_B(b_1)>t$. Hence
\begin{eqnarray*}
\mu_{A+B}([x_1, \ldots, x_n])&=&\mu_{A+B}([a_1+b_1, \ldots, x_n])\\
&\geq&\sup_{[x_1, \ldots, x_n]=[a, \ldots, x_n]+[b, \ldots, x_n]}\{\mu_A([a, \ldots, x_n])\wedge \mu_B([b, \ldots, x_n])\}\\
&\geq& \mu_A([a_1, \ldots, x_n])\wedge\mu_B([b_1, \ldots, x_n])\\
&\geq&( \mu_A(a_1)\vee \cdots \vee \mu_A(x_n))\wedge(\mu_B(b_1)\vee \cdots \vee \mu_B(x_n))\\
&>&t\\
&>&\mu_{A+B}([x_1, \ldots, x_n]).
\end{eqnarray*}
Contradiction. Thus $\mu_{A+B}([x_1, \ldots, x_n]\geq \mu_{A+B}(x_1)\vee \cdots \vee \mu_{A+B}(x_n)$. Finally, suppose $\lambda_{A+B}([x_1, \ldots, x_n])>\lambda_{A+B}(x_1)\wedge \cdots \wedge \lambda_{A+B}(x_n)$. Then there exists $s$ such that $\lambda_{A+B}([x_1, \ldots, x_n])>s>\lambda_{A+B}(x_1)\wedge \cdots \wedge \lambda_{A+B}(x_n)$. Without loss of generality we may assume $\lambda_{A+B}(x_1)=\lambda_{A+B}(x_1)\wedge \cdots \wedge \mu_{A+B}(x_n)$. Thus $\lambda_{A+B}([x_1, \ldots, x_n])>s>\inf_{x_1=a+b}\{\lambda_A(a)\vee\lambda_B(b)\}$, and so there exist $a_1,~b_1\in L$ such that $\lambda_{A+B}([x_1, \ldots, x_n])>s>\lambda_A(a_1)\vee\lambda_B(b_1)$. Therefore $\lambda_A(a_1)<s$ and $\lambda_B(b_1)<s$. Now
\begin{eqnarray*}
\lambda_{A+B}([x_1, \ldots, x_n])&=&\lambda_{A+B}([a_1+b_1, \ldots, x_n])\\
&\leq&\inf_{[x,~y]=[a, \ldots, x_n]+[b, \ldots, x_n]}\{\lambda_A([a, \ldots, x_n])\vee\lambda_B([b, \ldots, x_n])\}\\
&\leq& \lambda_A([a_1, \ldots, x_n])\vee\lambda_B([b_1,\ldots, x_n])\\
&\leq &(\lambda_A(a_1)\wedge \cdots \wedge \lambda_A(x_n))\vee(\lambda_B(b_1)\wedge \cdots \wedge \lambda_B(x_n))\\
&<&s\\
&<&\lambda_{A+B}([x_1, \ldots, x_n]).
\end{eqnarray*}
Contradiction. Thus $A+B$ is a complex intuitionistic fuzzy $n$-Lie ideal of $L$. \hfill$\Box$

In particular, if $A$ and $B$ are fuzzy Lie ideals of a Lie algebra $L$, then we obtain a special case of Shaqaqha's result \cite[Theorem 3.5]{shadi}.
\section{Direct Product of Intuitionistic fuzzy Lie $n$-subalgebras}
Let $A=(\mu_A, \lambda_A)$ and $B=(\mu_B, \lambda_B)$ be an intuitionistic fuzzy subsets of $L_1$ and $L_2$, respectively. Then the generalized cartesian product $A\times B$ is defined to be $A\times B=(\mu_A\times \mu_B, \lambda_A, \lambda_B)$ where
$$\mu_A\times \mu_B:L\times L\rightarrow [0,1];~(x, y)\mapsto \mu_A(x)\wedge \mu_B(y),$$
and
$$\lambda_A\times\lambda_B:L\times L\rightarrow [0, 1];~(x, y)\mapsto \lambda_A(x)\vee\lambda_B(y).$$
We have the following result.
\begin{Thm}\label{proifnls1}
	Let $A=(\mu_A, \lambda_A)$ and $B=(\mu_B, \lambda_B)$ be an intuitionistic fuzzy $n$-Lie subalgebras of $L_1$ and $L_2$, respectively.
	Then $A\times B$ is an intuitionistic fuzzy $n$-Lie subalgebra of $L_1\times L_2$.
\end{Thm}
{\it Proof.~}Note that $A\times B$ is an intuitionistic fuzzy subset of $L_1\times L_2$. Indeed if $x\in L_1$ and $y\in L_2$ such that $\mu_A(x)\leq \mu_B(y)$ and $\lambda_B(y)\geq \lambda_B(x)$, then 
$$(\mu_A\times\mu_B)(x, y)+(\lambda_A\times\lambda_B)(x, y) =\mu_A(x)+\lambda_B(y)\leq \mu_B(y)+\lambda_B(y)\leq 1.$$
Similarly, if $\mu_B(y)\leq \mu_A(x)$ and $\lambda_A(x)\geq \lambda_B(y)$. Let $(x_1, y_1), (x_2, y_2)\in L_1\times L_2$. Then the proofs for $(\mu_A\times\mu_B)((x_1, y_1)+ (x_2, y_2))\geq (\mu_A\times \mu_B)((x_1, y_1))\wedge  (\mu_A\times \mu_B)((x_2, y_2))$, $(\lambda_A\times\lambda_B)((x_1, y_1)+ (x_2, y_2))\leq \lambda_A\times \lambda_B((x_1, y_1))\vee  \lambda_A\times \lambda_B((x_2, y_2))$, $(\mu_A\times\mu_B)(c(x_1, y_1))\geq  (\mu_A\times \mu_B)((x_1, y_1))$, and $(\lambda_A\times\lambda_B)(c(x_1, y_1))\leq  (\lambda_A\times \lambda_B)((x_1, y_1))$ are similar to the proof of \cite[Lemma 2.1]{chen2}.\\
For $(x_1, y_1), (x_2, y_2), \ldots, (x_n, y_n)\in L\times L$, we have
\begin{eqnarray*}
	\mu_A([(x_1, y_1), (x_2, y_2), \ldots, (x_n, y_n)])&=& (\mu_A\times \mu_B)([x_1, x_2, \ldots, x_n], [y_1, y_2, \ldots, y_n])\\
	&=&\mu_A([x_1, x_2, \ldots, x_n])\wedge \mu_B([y_1, y_2, \ldots, y_n])\\
	&\geq& (\mu_A(x_1)\wedge\cdots \wedge \mu_A(x_n))\wedge(\mu_B(y_1)\wedge \cdots \wedge \mu_B(y_n))\\
	&=&(\mu_A\times \mu_B)((x_1, y_1))\wedge \cdots \wedge (\mu_A\times \mu_B)((x_n, y_2)),
\end{eqnarray*} 
and also
\begin{eqnarray*}
	\lambda_A([(x_1, y_1), (x_2, y_2), \ldots, (x_n, y_n)])&=& (\lambda_A\times \lambda_B)([x_1, x_2, \ldots, x_n], [y_1, y_2, \ldots, y_n])\\
	&=&\lambda_A([x_1, x_2, \ldots, x_n])\vee \lambda_B([y_1, y_2, \ldots, y_n])\\
	&\leq& (\lambda_A(x_1)\vee\cdots \vee \lambda_A(x_n))\vee(\lambda_B(y_1)\vee \cdots \vee \lambda_B(y_n))\\
	&=&(\lambda_A\times \lambda_B)((x_1, y_1))\vee \cdots \vee (\lambda_A\times \lambda_B)((x_n, y_2)).
\end{eqnarray*} 
\hfill $\Box$

However the direct product of two intuitionistic fuzzy Lie ideals of $n$-Lie algebras $L_1$ and $L_2$ is not nesaccary to be an intuitionistic fuzzy Lie ideal of the $n$-Lie algebra $L_1\times L_2$.
\section{On Lie Algebra Homomorphism of Intuitionistic Fuzzy $n$-Lie Algebras}
Let $L_1$ and $L_2$ be $n$-Lie algebras, $A=(\mu_A, \lambda_A)$ be an intuitionistic subset of $L_1$, and $f:L_1\rightarrow L_2$ be a function. Then the intuitionistic fuzzy subset $f(A)$ of $f(L_1)$ defined by $f(A)=(\mu_{f(A)}, \lambda_{f(A)})$ where $\mu_{f(A)}(y)=\sup_{x\in f^{-1}(y)}\{\mu_A(x)\}$ and $\lambda_{f(A)}(y)=\inf_{x\in f^{-1}(y)}\{\lambda_A(x)\}$ for $y\in f(L_1)$ is called the image of $A$ under $f$. Similarly, if $B=(\mu_B,\lambda_B)$ is an intuitionistic fuzzy subset of $L_2$, then the intuitionistic fuzzy set on $L_1$ is $f^{-1}(B)=(\mu_{f^{-1}(B)},\lambda_{f^{-1}(B)})$ where
$$\mu_{f^{-1}(B)}(x)=\mu_B(f(x)) ~\mathrm{and}~\lambda_{f^{-1}(B)}(x)=\lambda_B(f(x))$$
(see for example \cite{shadi}). The proofs of the following two results are omitted because they are routine and parallel to the corresponding results on intuitionistic fuzzy Lie algebras (\cite{Akram, shadi2}).
\begin{Thm}\label{homt1}
	Let $\varphi:L_1\rightarrow L_2$ be a $n$-Lie algebra homomorphism. If $B=(\mu_B,~\lambda_B)$ is an intuitionistic fuzzy $n$-Lie subalgebra (resp. ideal) on $L_2$, then the intuitionistic fuzzy set $\varphi^{-1}(B)$ is an intuitionistic fuzzy $n$-Lie subalgebra (rep. ideal) of $L_1$.
\end{Thm}
\begin{Thm}\label{homt2}
	Let $\varphi:L_1\rightarrow L_2$ be a $n$-Lie algebra homomorphism. If $A=(\mu_A,~\lambda_A)$ is an intuitionistic fuzzy $n$-Lie subalgebra (resp. ideal) on $L_1$, then the intuitionistic fuzzy set $\varphi(A)$ is an intuitionistic fuzzy $n$-Lie subalgebra (resp. ideal) of $\mathrm{im}(\varphi)$.	
\end{Thm}
\section{Intuitionistic Fuzzy Quotient $n$-Lie Algebras}
Let $L$ be an $n$-Lie algebra and $I$ be an ideal of $L$. Then the factor space $L/I=\{x+I~:~x\in I\}$ acquires an $n$-Lie algebra structure (called a quotient $n$-Lie algebra) by setting $$[x_1+I, x_2+I, \ldots, x_n+I]= [x_1,~ x_2, \ldots, x_n]+ I$$
for $x_1, x_2, \ldots, x_n\in L$.
In this paper we define and study the intuitionistic Fuzzy quotient $n$-Lie algebra by an intuitionistic fuzzy $n$-Lie ideal.
\begin{Def}\label{Q1}
Let $A=(\mu_A, \lambda_A)$ be an intuitionistic fuzzy $n$-Lie ideal of an $n$-Lie algebra $L$. Then for $x\in L$, the intuitionistic fuzzy subset $$x+A=(x+\mu_A, x+\lambda_A)$$
where
$$x+\mu_A:L\rightarrow[0,~1];~y\mapsto \mu_A(y-x)$$
and
$$x+\lambda_A:L\rightarrow[0,~1];~y\mapsto \lambda_A(y-x)$$
is called a coset (determined by $x, \mu_A$, and $\lambda_A$) of the intuitionistic fuzzy $n$-Lie ideal $A$.
\end{Def}
The case where $n=2$ was introduced and studied by Chen \cite{chen} in 2010. The set of all cosets of an intuitionistic fuzzy $n$-Lie ideal will be denoted by $L/A$. 
The following lemma proves that a coset may have many different
labels.
\begin{lemma}\label{Q2}
	Let $A=(\mu_A,~\lambda_A)$ be an intuitionistic fuzzy $n$-Lie ideal of an $n$-Lie algebra $L$, and let $x, y\in L$. The following statements are equivalent:
	\begin{itemize}
		\item[(i)] $x+A=y+A$,
		\item[(ii)] $\mu_A(x-y)=\mu_A(0)$ and $\lambda_A(x-y)=\lambda_A(0)$.
	\end{itemize}
\end{lemma}
{\it Proof.~}If $x,~y\in L$ with $x+A=y+A$, then $x+\mu_A=y+\mu_A$ and $x+\lambda_A=y+\lambda_A$. Consider $x+\mu_A=y+\mu_A$, then evaluating both sides for $x$ implies that $\mu_A(0)=\mu_A(x-y)$. Similarly $\lambda_A(x-y)=\lambda_A(0)$. Conversely, let  $\mu_A(x-y)=\mu_A(0)$ and $\lambda_A(x-y)=\lambda_A(0)$. Then for any $z\in L$, we have $(y+\mu_A)(z)=\mu_A(z-y)\geq\mu_A(z-x)\wedge\mu_A(x-y)=\mu_A(z-x)\wedge\mu_A(0)=\mu_A(z-x)=(x+\mu_A)(z)$.
Thus $y+\mu_A\geq x+\mu_A$. Also $\mu_A(y-x)=\mu_A(-(x-y))= \mu_A(x-y)=\mu_A(0)$. So we can prove that $x+\mu_A\geq y+\mu_A$ in the same way as above. Hence $x+\mu_A=y+\mu_A$. By almost the same aegument we can prove that $x+\lambda_A=y+\lambda_A$.\hfill$\Box$
\begin{Thm}\label{Q3}
	Let $A=(\mu_A,~\lambda_A)$ be an intuitionistic fuzzy $n$-Lie ideal of $L$ and $L/A$ be the set of all cosets of $L$ on $A$. Then the set $L/A$ is an $n$-Lie algebra under the following operations:
\begin{itemize}
	\item[(i)] $(x+A)+(y+A)=(x+y)+A$ for all $x, y\in L$,
	\item [(ii)] $\alpha(x+A)=\alpha x+A$ for all $x\in L$ and $\alpha \in F$,
	\item[(iii)] $[(x_1+A),(x_2+A), \ldots, (x_n+A)]=[x_1, x_2, \ldots, x_n]+A$ for all $x_1, x_2, \ldots, x_n\in L$.
\end{itemize}
\end{Thm}
{\it Proof.~}First, we show that the operations are well defined. Let $x,~y,~u$ and $v$ be elements in $L$ such that $x+A=y+A$ and $u+A=v+A$. Consequently $\mu_A(x+u-(y+v))=\mu_A((x-y)+(u-v))
\geq \mu_A(x-y)\wedge\mu_A(u-v)
=\mu_A(0)$ (because $\mu_A(x-y)=\mu_A(0)~\mathrm{and}~\mu_A(u-v)=\mu_A(0)$).
As $\mu_A(0)\geq \mu_A(x+u-(y+v))$, we have $\mu_A(x+u-y-v)=\mu_A(0)$. Almost the same argument one can obtain that $\lambda_A(x+u-y-v)=\lambda_A(0)$ Therefore $(x+u)+A= (y+v)+A$. Also, $\mu_A(\alpha x-\alpha y)
\geq \mu_A( x-y)=\mu_A(0)$ and $\lambda_A(\alpha x-\alpha y)
\leq \lambda_A( x-y)=\lambda_A(0)$. Hence $\alpha(x+A)=\alpha(y+A)$. If $x_1, x_2, \ldots, x_n\in L$ such that $x_i+A= y_i+A$ ($i=1, \ldots, n$), then
\begin{eqnarray*}
\mu_A([x_1,x_2, \ldots, x_n]-[y_1, y_2, \ldots, y_n])%&=&\mu_A([x_1, x_2, \ldots, x_{n-1}, x_n]-[y_1, x_2, \ldots, x_n]\\
%&&+[y_1, x_2, \ldots, x_n]-[y_1, \ldots, y_{n-1}, y_n])\\
%&=&\mu_A([x_1, \ldots, x_{n-1}, x_n]+[-y_1, x_2, \ldots, x_n]\\
%&&+[y_1, x_2, \ldots, x_n]-[y_1, y_2, \ldots, y_n])\\
&=&\mu_A([x_1-y_1,x_2, \ldots, x_n]\\
&&+[y_1, x_2, \ldots, x_n]-[y_1, y_2, \ldots, y_n])\\
&\geq&\mu_A([x_1-y_1,x_2, \ldots, x_n])\\
&&\wedge \mu_A([y_1, x_2, \ldots, x_n]-[y_1, y_2, \ldots, y_n])\\
%&\geq&(\mu_A(x_1-y_1)\vee \mu_A(x_2)\vee \cdots \vee \mu_A(x_n))\\
%&&\wedge \mu_A([y_1, x_2, \ldots, x_n]-[y_1, \ldots, y_{n-1}, y_n])\\
&\geq&\mu_A(0)\wedge \mu_A([y_1, x_2, \ldots, x_n]-[y_1, y_2, \ldots, y_n])
\end{eqnarray*}
Now,
\begin{eqnarray*}
 \mu_A([y_1, x_2, \ldots, x_n]-[y_1, y_2, \ldots, y_n])%&=& \mu_A([y_1, x_2, \ldots, x_n]+[y_1, -y_2, x_3, \ldots, x_n]\\
 %&&+[y_1, y_2, x_3, \ldots, x_2]-[y_1, y_2, y_3, \ldots, y_n])\\
 &=& \mu_A([y_1, x_2-y_2, x_3, \ldots, x_n]\\
 &&+[y_1, y_2, x_3, \ldots, x_2]-[y_1, y_2, y_3, \ldots, y_n])\\
%  &\geq& \mu_A([y_1, x_2-y_2, x_3, \ldots, x_n])\\
 %&&\wedge\mu_A([y_1, y_2, x_3, \ldots, x_2]-[y_1, y_2, y_3, \ldots, y_n])\\
 &\geq& \mu_A(y_1)_\vee\mu_A(x_2-y_2)\vee \mu_A(x_3)\vee \cdots\mu_A(x_n)\\
 &&\wedge\mu_A([y_1, y_2, x_3, \ldots, x_2]-[y_1, y_2, y_3, \ldots, y_n])\\
 &\geq&\mu_A(0)\wedge\mu_A([y_1, y_2, x_3, \ldots, x_2]-[y_1, y_2, y_3, \ldots, y_n]).
\end{eqnarray*}
So 
$$\mu_A([x_1,x_2, \ldots, x_n]-[y_1, y_2, \ldots, y_n])\geq \mu_A(0)\wedge\mu_A(0) \wedge\mu_A([y_1, y_2, x_3, \ldots, x_2]-[y_1, y_2, y_3, \ldots, y_n]).$$
Continuing in the same way, we obtain 
$$\mu_A([x_1,x_2, \ldots, x_n]-[y_1, y_2, \ldots, y_n])\geq \mu_A(0)\wedge\mu_A(0) \wedge \cdots \wedge \mu_A(0) ~(n~\mathrm{times}).$$ Hence $\mu_A([x_1,x_2, \ldots, x_n]-[y_1, y_2, \ldots, y_n])=\mu_A(0)$. Using the same method one can prove that $\lambda_A([x_1, x_2, \ldots, x_n]-[y_1, y_2, \ldots, y_n])=\lambda_A(0)$. Therefore  $[x_1, x_2, \ldots, x_n]+A=[y_1, y_2, \ldots, y_n]+A$. Now it is straightforward to see that the product on $L/A$ is an $n$-linear map satisfying the Filippov identity.\hfill$\Box$

The $n$-Lie algebra constructed in Theorem \ref{Q3} is called intuitionistic fuzzy quotient $n$-Lie algebra of $L$ by $A$.


\begin{thebibliography}{99}
	\bibitem{Abdukhalikov} K.S. Abdukhalikov, M.S. Tulenbaev, U.U. Umirbaev, On fuzzy subalgebras, \textit{Fuzzy Sets and Systems}, {\bf 93} (1998), 257-262.
	\bibitem{Akram book} M. Akram, \textit{Fuzzy Lie algebras},	Springer Nature Singapore Pte Ltd. (2018).	
	\bibitem{Akram} M. Akram, Intuitionistic fuzzy Lie subalgebras, {\it Southeast Asian Bulletin of Mathematics} {\bf 31} (2007), 843-855.
	\bibitem{Akram1} M. Akram, Intuitionistic (S,T)-fuzzy Lie ideals of Lie algebras, {\it Quasigroups Relat. Systems} {\bf 15}
(2007), 201-215.
    \bibitem{Akram2} M. Akram, Intuitionistic fuzzy Lie ideals of Lie algebras, {\it Int. Journal of Fuzzy Math.}, {\bf 6} (4)
(2008), 991-1008.
    \bibitem{Akram3} M. Akram, Co-fuzzy Lie superalgebras over a co-fuzzy field, {\it World Applied Sciences Journal}, {\bf 7}
(2009), 25-32.
	\bibitem{Alekseevsky} D. Alekseevsky and P. Guha, On decomposability of Nambu-Poisson Tensor, {\it Acta
Mathematica Universitatis Comenianae}, {\bf 65} (1996), 1-9
%	\bibitem{Alkouri} A.S. Alkouri and A. Salleh, Complex intuitionistic fuzzy sets, \textit{in Proceedings of the International Conference on Fundamental and Applied Sciences}  (ICFAS ’12), {\bf1482} (2012), 464-470.
%	\bibitem{Salleh} A.S. Alkouri and A. Salleh, Some operations on complex Atanassov’s intuitionistic fuzzy sets, \textit{AIP Conference Proceedings}, {\bf 1571} (2013), 987.
	%\bibitem{Alsarahead} M. Alsarahead and A. Ahmad, Complex fuzzy subgroups, \textit{Applied Mathematical Sciences}, {\bf41}(11) (2017), 2011--2021.
	%\bibitem{Alsarahead2} M. Alsarahead and A. Ahmad, Complex fuzzy subrings, \textit{International Journal of Pure	and Applied Mathematics}, {\bf117}(4) (2017), 563--573.
	%\bibitem{Alsarahead3} M. Alsarahead and A. Ahmad, Complex intuitionistic fuzzy subgroups, \textit{submitted to	Italian Journal of Pure and Applied Mathematics} (2017).
	%\bibitem{Alsarahead4} M. Alsarahead and A. Ahmad, Complex intuitionistic fuzzy subrings, \textit{Borneo Science} {\bf38}(2) (2018), 24-37.
	\bibitem{ATANASSOV} K.T Atanassov, Intuitionistic fuzzy sets, {\it Fuzzy Sets and Systems} {\bf 20} (1986), 87-96.
	\bibitem{Bahturin 1} Y. Bahturin, \textit{Identical relations in Lie algebras}, VNU Science Press, b.v., Utrecht, 1987.
	\bibitem{Biswas} R. Biswas, Intuitionistic fuzzy subgroups, {\it Math. Forum} {\bf 10} (1989),  37-46.
	\bibitem{chen2} W. Chen, S. Zhang, Intuitionistic fuzzy Lie sub-superalgebra and intuitionistic fuzzy ideals, {\it Computer and Mathematics with Applications} {\bf 58} (2009), 1645-1661.
	\bibitem{chen} W. Chen, Intuitionistic fuzzy quotient Lie superalgebras, {\it International Journal of Fuzzy Systems},  {\bf12} (4), (2010), 330-339.
	\bibitem{Davvaz} B. Davvaz and WA. Dudek, Fuzzy $n$-Lie algebras, {\it J Generalized Lie Theory Appl},{\bf 11} (2017), 1-6.
	\bibitem{Dudek} W. A. Dudek, Fuzzifications of $n$-ary groupoids, {\it Quasigroups and Related Systems},{\bf 7} (2000), 45-66.
	\bibitem{Filippov} Filippov, \textit{$n$-ary Lie algebras}, (Russian) Sibirsk Mat Zh, {\bf 26}(6) (1985), 126--140.
	\bibitem{Gautheron} P. Gautheron, Simple facts concerning Nambu algebras, {\it Commun. Math. Phys.}, {\bf 195} (1998), 417-34
	\bibitem{Ho} P. Ho, M. Chebotar, and W. Ke, On skew-symmetric maps on Lie algebras, {\it Proc. Royal
Soc. Edinburgh}, {\bf 133A} 2003, 1273-1281.
	\bibitem{Karin} K. Erdmann and M.J. Wildon, \textit{Introduction to Lie algebra}, Springer Undergraduate Mathematics Series. Spinger-Verlag London Limited (2006).
	%\bibitem{Joseph} J.A. Gallian, \textit{Contemporary Abstract Algebra}, Eighth Edition, Boston, (2012).
	\bibitem{Jacobson} N. Jacobson, \textit{Lie Algebras}, Wiley, New York, (1962). 
	\bibitem{Katsaras} A.K. Katsaras and D.B. Liu, Fuzzy vector spaces and fuzzy topological vector spaces, \textit{J. Math. Anal. Appl.}, {\bf 58} (1997), 135-146. 
	\bibitem{Kim} C. Kim and D. Lee, Fuzzy Lie ideals and fuzzy Lie subalgebras, {\it Fuzzy Sets and Systems} {\bf 94} (1998), 101-107.
	\bibitem{Kondo} M. Kondo and W. A. Dudek, On the transfer principal in fuzzy theorey, {\it Mathware Soft Computing} {\bf 12} (2005), 41-55.
	\bibitem{Marmo} G. Marmo, G. Vilasi and A. M. Vinogradov, The local structure of $n$-Poisson and $n$-Jacobi manifolds, {\it J. Geom. Phys.}, {\bf 25} (1998), 141-82
	\bibitem{Michor} P. W. Michor and A. M. Vinogradov, $n$-ary and associative algebras, {\it Rend. Sem. Mat. Univ. Pol. Torino}, {\bf 53} (1996), 373-92.
	\bibitem{Nakanishi} N. Nakanishi, On Nambu-Poisson manifolds {\it Rev. Math. Phys.}, {\bf 10} (1998), 499-510.
	\bibitem{Nambu} Y. Nambu, Generalized Hamiltonian Dynamics, {\it Physics. Rev.}, {\bf D7} (1973), 2405-2412.
	%\bibitem{Schaum} S. Lipschutz and M.L Lipson, \textit{Schaum's outline of Theory and Problems of Linear Algebra}, Fourth Edition, McGraw-Hill, (2009).
	\bibitem{maik} D.S. Mailk and J.N. Mordeson, Fuzzy vector spaces, {\it lnformation Sciences} {\bf 55} (1991), 271-281.
	\bibitem{Papadopoulos} G, Papadopoulos, M2-branes, $3$-Lie algebras and Plucker relations, arXiv: 0804.
2662[hep-th]
	\bibitem{Takhtajan} L. Takhtajan, On foundation of the generalized Nambu mechanics, {\it Commun. Math.
Phys.}, {\bf 160} (1993), 295-315.
	\bibitem{Ramot} D. Ramot, M. Friedman, G. Langholz and A. Kandel, Complex fuzzy sets, \textit{IEEE Transaction on Fuzzy Systems}, {\bf 10}(2) (2002), 171-186.
	\bibitem{Rosenfeld} A. Rozenfeld, Fuzzy groups, \textit{J. Math. Anal. Appl.}, {\bf 35} (1971), 512–517.
	\bibitem{shadi} S. Shaqaqha, Complex fuzzy Lie algebras, {\it Jordan Journal of Mathematics and Statistics}, {\bf 13}(2) (2020), 227-244.
	\bibitem{shadi2} S. Shaqaqha and M. Al-Deiakeh, Complex intuitionistic fuzzy Lie subalgebras, preprint.
	\bibitem{Vinogradov} A. Vinogradov and M. Vinogradov, On multiple generalizations of lie algebras and
poisson manifolds,{\it American Mathematical Society, Contemp. Math.}, {\bf 219} (1998), 273-287.
	\bibitem{Yehia1} S.E. Yehia,  Fuzzy ideals and fuzzy subalgebras of Lie algebra, {\it Fuzzy Sets and Systems} {\bf 80} (1996), 237-244.
    \bibitem{Yehia2} S.E. Yehia, The adjoint representation of fuzzy Lie algebras, {\it Fuzzy Sets and Systems} {\bf 119} (2001), 409-417.
	\bibitem{Zadeh} L. Zadeh, Fuzzy sets, {\it Inform. Control} {\bf 8} (1965), 338-358.
\end{thebibliography}
\end{document}